 \documentclass[envcountsect, envcountsame]{svmult}
 
\usepackage{mathptmx}       % selects Times Roman as basic font
\usepackage{helvet}         % selects Helvetica as sans-serif font
\usepackage{courier}        % selects Courier as ty1pewriter font
\usepackage{makeidx}         % allows index generation
\usepackage{graphicx}        % standard LaTeX graphics tool
\usepackage{multicol}        % used for the two-column index
\usepackage[bottom]{footmisc}% places footnotes at page bottom
\makeindex             % used for the subject index
\usepackage{amstext,amsfonts,amssymb,amscd,amsbsy,amsmath}
\usepackage{tikz-cd}	
\usepackage{tikz}

\newcommand{\A}{\mathbb{A}}

\newcommand{\NN}{\mathbb{N}}
\newcommand{\CC}{\mathbb{C}}
\newcommand{\PP}{\mathbb{P}}

\newcommand{\Hrk}{\mathrm{Hrk}}

\smartqed

\begin{document}

\title*{Minkowski sums and Hadamard products of algebraic varieties}
\author{Netanel Friedenberg, Alessandro Oneto, and Robert L. Williams}
\institute{
Netanel Friedenberg \at Yale University,\\ 10 Hillhouse Ave-Ste 442\\ PO Box 208283\\ New Haven, CT  06520-8283, United States\\ \email{netanel.friedenberg@yale.edu} \and 
Alessandro Oneto \at INRIA Sophia Antipolis M\'editerran\'ee, \\ 2004 Route de Lucioles, \\ 06902 Sophia Antipolis, France, \\ \email{alessandro.oneto@inria.fr} \and 
Robert L. Williams \at Texas A\&M University, Department of Mathematics \\
Mailstop 3368 \\
College Station, TX 77843-3368 United States \\ \email{rwilliams@math.tamu.edu}}
\maketitle

\abstract*{We study Minkowski sums and Hadamard products of algebraic varieties. Specifically we explore when these are varieties and examine their properties in terms of those of the original varieties. This project was inspired by Problem~5 on Surfaces in \cite{Sturmfels}.}

\abstract{We study Minkowski sums and Hadamard products of algebraic varieties. Specifically we explore when these are varieties and examine their properties in terms of those of the original varieties. This project was inspired by Problem~5 on Surfaces in \cite{Sturmfels}.
}

\section{Introduction}

In algebraic geometry we have several constructions to build new algebraic varieties from given ones. %A very classical and well-studied example is the {\it join} of two varieties, i.e., the closure of the union of lines connecting the varieties, as well as secant varieties, rational normal scrolls, and Segre products. 
Examples of classical, well-studied constructions are joins, secant varieties, rational normal scrolls, and Segre products. 
In these cases, it is very interesting to understand geometric properties, e.g., the dimension and the degree, of the variety constructed in terms of those of the original varieties. In this chapter we focus on the {\it Minkowski sum} and the {\it Hadamard product} of algebraic varieties. These are constructed by considering the entry-wise sum and multiplication, respectively, of points on the varieties. Due to the nature of these operations, there is a remarkable difference between the affine and the projective case.

The entry-wise sum is not well-defined over projective spaces. For this reason, we consider only Minkowski sums of affine varieties. However, in the case of affine cones, the Minkowski sum corresponds to the classical join of the corresponding projective varieties. Conversely, we focus on Hadamard products of projective varieties and, in particular, of varieties of matrices with fixed rank. This is because these Hadamard products parametrize interesting problems related to algebraic statistics and quantum information.
 
\smallskip
 
Our original motivating question was the following.

\begin{question}{\it Which properties do the Minkowski sum and the Hadamard product have with respect to the properties of the original varieties? In particular, what are their dimensions and degrees?}
\end{question}

\noindent We now introduce these constructions. 
 We work over an algebraically closed field $\Bbbk$. We will add extra assumptions on $\Bbbk$ when needed. We use the notation $\Bbbk^\times:=\Bbbk\setminus\{0\}$.

\begin{definition}%[affine case]
 Let $X,Y \subset \A^n$ be affine varieties. We define the \emph{Minkowski sum} of $X$ and $Y$, denoted $X+Y$, as the Zariski closure of the image of $X\times Y$ under the entry-wise summation map 
\begin{eqnarray*}
  \phi_+: & \A^n \times \A^n & \rightarrow ~~~~~~~~~~~~~~ \A^n, \\
  & ((a_1,\ldots,a_n),(b_1,\ldots,b_n)) &\mapsto (a_1+b_1,\ldots,a_n+b_n)
\end{eqnarray*}
 \end{definition}

Note that taking the Zariski closure of $\phi_+(X\times Y)$ is necessary to construct an algebraic variety, as explained in Example \ref{example:TwoHyperbolas}.
 
 \medskip
 
\noindent As far as we know there is no literature about Minkowski sums of varieties. We compute the dimension and degree of Minkowski sums of generic affine varieties.
\setcounter{section}{3}
\setcounter{theorem}{8}
\begin{theorem}
 Let $X,Y \subset \A^n$ be varieties. Then, for $X$ and $Y$ in general position, $\dim(X+Y) = \min\{\dim(X)+\dim(Y),n\}$.
\end{theorem}
\setcounter{theorem}{11}
\begin{corollary}
Suppose $\Bbbk$ has characteristic other than $2$. Let $X,Y\subset \A^n$ be varieties whose projective closures $\overline{X},\overline{Y}\subset\PP^n$ are contained in complementary linear subspaces; equivalently, $X,Y$ are contained in disjoint affine subspaces which are not parallel. Then for generic $\alpha\in\Bbbk^\times$, $\deg(\alpha X+Y) = \deg(X)\deg(Y)$.
\end{corollary}
A crucial observation in our computations is that the Minkowski sum of affine varieties disjoint at infinity can be described in terms of the join of their projectivizations, see Proposition \ref{prop:CayleyTrick} and Remark \ref{remark:Cayley}. This is a construction inspired by the combinatorial Cayley trick used to construct Minkowski sums of polytopes. 

\setcounter{section}{1}
\setcounter{theorem}{2}
 
%\noindent {\sc Hadamard product of projective varieties.}
%\medskip

%\begin{definition}[affine case]
% Let $X,Y \subset \A^n$ be affine varieties. We define the {\it Hadamard product} of $X$ and $Y$, denoted by $X\star Y$, as the Zariski closure of the image of $X\times Y$ under the entry-wise multiplication map 
%\begin{eqnarray*}
%  \phi_\star: & \A^n \times \A^n & \rightarrow ~~~~~~~~~~ \A^n, \\
%  & ((a_1,\ldots,a_n),(b_1,\ldots,b_n)) &\mapsto (a_1b_1,\ldots,a_nb_n)
%\end{eqnarray*}
% \end{definition}
% For Hadamard products, we can directly extend the definition to projective varieties. 

%\setcounter{definition}{1}
\begin{definition}
  Let $X,Y \subset \PP^n$ be projective varieties. We define the {\it Hadamard product} of $X$ and $Y$, denoted by $X\star Y$, as the Zariski closure of the image of $X\times Y$ under the map 
\begin{eqnarray*}
  \phi_\star: & \PP^n \times \PP^n & \dashrightarrow ~~~~~~~~~~~~~~~~~~~~ \PP^n, \\
  & ([a_0:\ldots:a_n],[b_0:\ldots:b_n]) &\mapsto [a_0b_0:a_1b_1:\ldots:a_nb_n].
\end{eqnarray*}
\end{definition} 
Let $\{x_0,\ldots,x_n\}$ be the homogeneous coordinates over $\PP^n$. The map $\phi_\star$ is not defined over the union of coordinate spaces $H_{I} \times H_{I^c}$, where $I \subset \{0,\ldots,n\}$, $I^c$ is its complement, and $H_I$ is the linear space defined by $\{x_i = 0 ~|~ i \in I\}$.
%In this case, the locus of non-definition of the map $\phi_\star$ is equal to the pairs of points 
%such that $a_ib_i = 0$, for all $i = 0,\ldots,n$. We call them {\it mutually exclusive points}.

Thus, the Hadamard product of projective varieties $X,Y \subset \PP^n$ is
$$
 X\star Y := \overline{\{p\star q ~:~ p\in X,~q\in Y,~p\star q \text{ is defined}\}} \subset \PP^n,
$$
where $p\star q := [p_0q_0:\ldots:p_nq_n]$ is the point obtained by entry-wise multiplication of the points $p = [p_0:\ldots:p_n]$ and $q = [q_0:\ldots:q_n]$. Also in this construction the operation of closure is crucial, as we show in Example \ref{example:Hadamard_closure}.

\medskip
\noindent %Hadamard products of projective varieties arose recently in \cite{CMS,PHYS}, \textcolor{blue}{in relation to some particular model in Statistics and Quantum Information.} 
In \cite{BCK}, the authors studied the geometry of Hadamard products, with a particular focus on the case of linear spaces. This work has been continued in \cite{BCFL}. 

In particular, we are interested in studying Hadamard products of varieties of matrices. The Hadamard product of matrices is a classical operation in matrix analysis \cite{H}. Its most relevant property is that it is closed on positive matrices. The Hadamard product of tensors appeared more recently in quantum information \cite{PHYS} and in statistics \cite{CMS,MM16}. In the latter, the authors studied restricted Boltzmann machines which are statistical models for binary random variables where some are hidden. From a geometric point of view, this reduces to studying {\it Hadamard powers} of the first secant variety of Segre products of copies of $\PP^1$. An interesting question is to understand how to express matrices as Hadamard products of small rank matrices. We call these expressions {\it Hadamard decomposition}. We define {\it Hadamard ranks} of matrices by using a multiplicative version of the usual definitions used for additive tensor decompositions. The study of Hadamard ranks is related to the study of Hadamard powers of secant varieties of Segre products of projective spaces.

In Section \ref{sec:Hadamard}, we focus in particular on the dimension of these Hadamard powers. We define the expected dimension and, consequently, we define the expected $r$-th Hadamard generic rank, i.e., the expected number of rank $r$ matrices needed to decompose the generic matrix of size $m\times n$ as their Hadamard product. It is 
$$
 {\rm exp}.\Hrk^\circ_r(m,n) = \left\lceil\frac{\dim\PP({\rm Mat}_{m,n}) - \dim(X_1)}{\dim(X_r) - \dim(X_1)}\right\rceil = \left\lceil\frac{mn-(m+n-1)}{r(m+n-r)-m-n+1}\right\rceil.
$$
We confirm this is correct for square matrices of small size using \textit{Macaulay2}.

% Our original motivating question was the following.
%\begin{question}{\it Which properties do the Minkowski sum and the Hadamard product have with respect to the properties of the original varieties? In particular, what are their dimensions and degrees?}
%\end{question}

%\noindent Before exploring these operations further, we make the following observation:

\bigskip
The paper is structured as follows. In Section \ref{sec:experiments}, we present some explicit computations of these varieties. We use both \emph{Macaulay2} \cite{M2} and \emph{Sage} \cite{sage}. These computations allowed us to conjecture some geometric properties of Minkowski sums and Hadamard products of algebraic varieties. 
In Section \ref{sec:Minkowski}, we analyze Minkowski sums of affine varieties. %In particular, we prove, under genericity conditions, that the dimension of the Minkowski sum is the sum of the dimensions and that its degree is the product of the degrees.
In particular, we prove that, under genericity conditions, the dimension of the Minkowski sum is the sum of the dimensions and we investigate the degree of the Minkowski sum.
 In Section \ref{sec:Hadamard}, we study Hadamard products and Hadamard powers of projective varieties. In particular, we focus on the case of Hadamard powers of projective varieties of matrices of given rank. We introduce the notion of Hadamard decomposition and Hadamard rank of a matrix. These concepts may be viewed as the multiplicative versions of the well-studied {\it additive decomposition} of tensors and {\it tensor ranks}.
%
%\section{Existing literature}
%
%From the literature, it seems that not much is known.
%
%\paragraph{\sc Minkowski sum.} We couldn't find anything about geometric properties of Minkowski sums.
%
%\paragraph{\sc Hadamard product.} In \cite{CMS}, Cueto, Morton and Sturmfels introduced Hadamard products in relation to some machine learning model (that we do not know anything about). In \cite{BCK}, Bocci, Carlini and Kileel started to investigate deeper the Hadamard product of projective varieties with particular focus on linear spaces.
\section{Experiments}\label{sec:experiments}
Problem~5 on Surfaces in \cite{Sturmfels} asked the following:
 \begin{center}
  {\it Compute the Minkowski sum and the Hadamard product \\  
  of two random circles in $\mathbb{R}^3$. Try other curves.}
 \end{center} 
 In order to compute Minkowski sums and Hadamard products of circles and other  curves, we used the algebra softwares \emph{Macaulay2} and \emph{Sage} to obtain equations and nice graphics. These also aided our general understanding of the geometric properties of these constructions.
Via elimination theory, we can compute the ideals of Minkowski sums and Hadamard products. This is the script in \emph{Macaulay2} to do so.
\begin{small}
\begin{verbatim}
 R = QQ[z_1..z_n,
        x_1..x_n,y_1..y_n];
 I = ideal( ... ); -- ideal of X in variables x_i;
 J = ideal( ... ); -- ideal of Y in variables y_i;
 
---- construct the ideals of graphs of the maps 
---- phi_+ and phi_star
 S = I + J + ideal(z_1-x_1-y_1,...,z_n-x_n-y_n);
 P = I + J + ideal(z_1-x_1*y_1,...,z_n-x_n*y_n);
 Msum = eliminate(toList{x_1..x_n | y_1..y_n}, S);
 Hprod = eliminate(toList{x_1..x_n | y_1..y_n}, P);
\end{verbatim}
\end{small}

With \emph{Sage}, we produced graphics of the real parts of Minkowski sums and Hadamard products of curves in $\A^3$. This is the script we used.

\begin{figure}
  \begin{center}
    \includegraphics[scale=0.4]{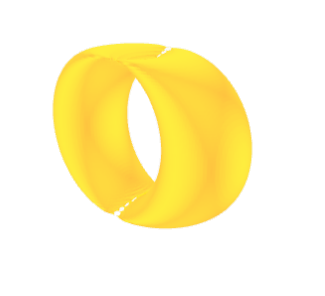}
    \caption{Minkowski sum of a circle of radius $1$ in the $x,y$-plane and a circle of radius $2$ in the $x,z$-plane. This is a degree four surface.} \label{sumCircles}
  \end{center}
\end{figure}

\begin{figure}
  \begin{center}
    \includegraphics[scale=0.4]{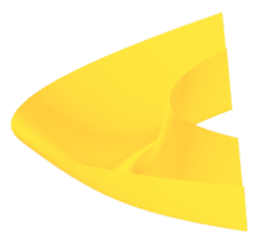}\\%($z^4 - 2yz^2 + y^2 - x=0$)
    \caption{Minkowski sum of the two parabolas $x=y^2$ in the $x,y$-plane and $y=z^2$ in the $y,z$-plane. This is a degree four surface.}\label{sumParabolas}
  \end{center}
\end{figure}

\begin{figure}
  \begin{center}
    \includegraphics[scale=0.45]{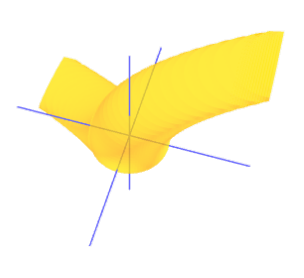}%($z^6 - 2xz^3 + z^4 - 2yz^2 + x^2 + y^2 - 1=0$)
    \includegraphics[scale=0.4]{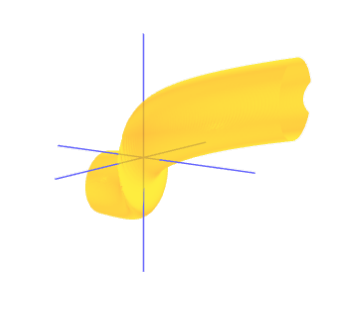}
    \includegraphics[scale=0.33]{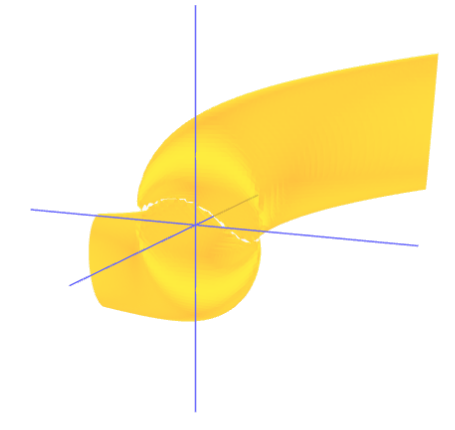}\\
    \caption{Minkowski sum of the twisted cubic with the unit circle in the $x,y$-plane, $y,z$-plane, and $x,z$-plane from left to right, respectively. Each of these are a degree six surface.} \label{sumTwistedCubicCircle}
  \end{center}
\end{figure}

\begin{figure} 
  \begin{center}
    \includegraphics[scale=0.4]{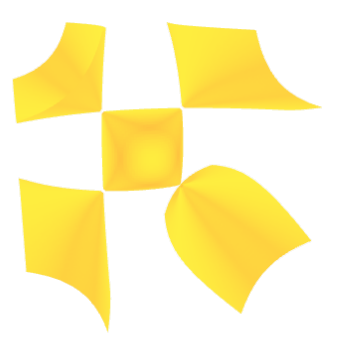}
    \includegraphics[scale=0.4]{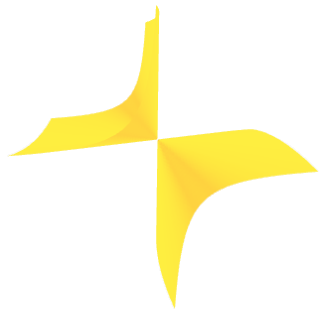}\\
    \caption{The real part of the Hadamard product of (left) the unit circle in the $z=1$ plane and the unit circle in the $y=1$ plane and (right) the circles $x^2+(y+z)^2=1,z-y=1$ and $x^2+(z-y)^2=1,y+z=1$. The surface on the left is of degree four, and the surface on the right is of degree two.}\label{hadamardProds}
  \end{center}
\end{figure}

\begin{small}
\begin{verbatim}
 A.<x1,x2,x3,y1,y2,y3,z1,z2,z3>=QQ[]
 I=( ... )*A # ideal of X in the variables x;
 J=( ... )*A # ideal of Y in the variables y;

# construct the ideals defining the graphs of the maps
# phi_+ and phi_star
 S = I + J + (z1-(x1+y1),z2-(x2+y2),z3-(x3+y3))*A
 P = I + J + (z1-(x1*y1),z2-(x2*y2),z3-(x3*y3))*A
 
 MSum = S.elimination_ideal([x1,x2,x3,y1,y2,y3])
 HProd = P.elimination_ideal([x1,x2,x3,y1,y2,y3])

# Assuming we get a surface, take the one generator of 
# each ideal.
 MSumGen=MSum.gens()[0]
 HProdGen=HProd.gens()[0]

# We plot these surfaces. 
# Because MSumGen and HProdGen are considered as elements of A, 
# which has 9 variables, they take 9 arguments.
 var('z1,z2,z3')
 implicit_plot3d(MSumGen(0,0,0,0,0,0,z1,z2,z3)==0, 
                (z1, -3, 3), (z2, -3,3), (z3, -3,3))
 implicit_plot3d(HProdGen(0,0,0,0,0,0,z1,z2,z3)==0, 
                (z1, -3, 3), (z2, -3,3), (z3, -3,3))
\end{verbatim}
\end{small}
In Figures \ref{sumCircles}, \ref{sumParabolas}, \ref{sumTwistedCubicCircle} and \ref{hadamardProds} are some of the pictures we obtained. These experiments gave us a first idea about the properties of Minkowski sums and Hadamard products.
\begin{note}\label{note:SumAndProductBasicProperties}
The fact that $X+Y$ and $X\star Y$ are the closures of images of $X\times Y$ under the maps $\phi_+$ and $\phi_\star$ immediately gives us that 
\begin{center}
$\dim(X+Y),\dim(X\star Y)\leq \dim(X)+\dim(Y)$
\end{center} and that if $X$ and $Y$ are irreducible, so are $X+Y$ and $X\star Y$.
\end{note}
Also, the fact that $X+Y$ and $X\star Y$ are linear projections of $X\times Y\subset\A^n\times\A^n$ and $X\times Y\subset\PP^n\times\PP^n\subset\PP^{n^2+2n}$, respectively, leads us to expect other geometric properties of $X+Y$ and $X\star Y$. 
%We immediately get that $\dim(X+Y),\dim(X\star Y)\leq\dim (X)+\dim(Y)$. In fact, 
%\textcolor{blue}{Since a linear projection from a linear space $L$ is generically one-to-one over a variety of at least complimentary dimension in generic position,}
Because the projection of a variety $Z\subset\PP^N$ in generic position from a linear space $L$ with $\dim(Z)+\dim(L)<N-1$ is generically one-to-one,
 we na\"ively expect that, for $X$ and $Y$ in general position, with $\dim(X)+\dim(Y)<n$,
\begin{center}$\dim(X+Y),\dim(X\star Y)=\dim(X)+\dim(Y),$

\smallskip
$\deg(X+Y)=\deg(X\times Y)=\deg(X)\deg(Y),~~~(X\times Y \subset \A^n)$

\smallskip
$\deg(X\star Y)=\deg(X\times Y)= \binom{\dim(X)+\dim(Y)}{\dim(X)}\deg(X)\deg(Y),~(X\times Y\subset\PP^{n^2+2n}).$
\end{center}
These expectations, however, do not follow directly from the projections of the varieties in general position because, even for $X$ and $Y$ in general position, $X\times Y$ is not in general position. Hence, we need further analysis as in the following sections.

\section{Minkowski sum of affine varieties}\label{sec:Minkowski}
Recall the definition of the Minkowski sum of two affine varieties $X$ and $Y$ as the closure of the image of $X\times Y$ under the map
\begin{eqnarray*}
  \phi_+: & \A^n \times \A^n & \rightarrow ~~~~~~~~~~~~~~ \A^n, \\
  & ((a_1,\ldots,a_n),(b_1,\ldots,b_n)) &\mapsto (a_1+b_1,\ldots,a_n+b_n)
\end{eqnarray*}
The operation of closure is needed in order to get an algebraic variety. Indeed, we can give an example where $X+_{\mathrm{set}}Y:=\phi_+(X\times Y)=\{p+q \mid p\in X,q\in Y\}$, the {\it setwise Minkowski sum} of $X$ and $Y$, is not closed.

\begin{example}\label{example:TwoHyperbolas}
In the affine plane $\A^2$ with coordinates $\{x,y\}$,
consider the plane curves $X=\{xy=1\}$ and $Y=\{xy=-1\}$. 
We claim that $\phi_+(X\times Y)$ contains the torus $(\Bbbk^\times)^2$, so $X+Y=\A^2$. For if $(\alpha,\beta)\in (\Bbbk^\times)^2$ and $(p,q)\in X\times Y$ then $\phi_+(p,q)=(\alpha,\beta)$ if and only if we can write $p$ and $q$ in the forms $p=(\alpha+a,\frac{1}{\alpha+a})$ and $q=(-a,\frac{1}{a})$ for some scalar $a\neq0,-\alpha$ such that $\beta=\frac{1}{\alpha+a}+\frac{1}{a}$. Clearing denominators in this last expression, we find the requirement is that $\beta a(\alpha+a)=a+(\alpha+a)$, i.e.\ $a$ is a zero of the quadratic polynomial $f_{\alpha,\beta}(t)=\beta t^2+(\beta\alpha-2)t-\alpha$. Note that $f_{\alpha,\beta}(0)=-\alpha$ and $f_{\alpha,\beta}(-\alpha)=\beta\alpha^2-\beta\alpha^2+2\alpha-\alpha=\alpha$. So if we let $a$ be a zero of $f_{\alpha,\beta}$ then with $p$ and $q$ as above we have $\phi_+(p,q)=(\alpha,\beta)$.

On the other hand $X+_{\mathrm{set}}Y$ is not all of $\A^2$. If the characteristic of the base field is not 2, then $X+_{\mathrm{set}}Y$does not contain the origin (though it does contain the punctured axes). In characteristic 2 we find that $X+_{\mathrm{set}}Y$ contains no point of the punctured axes $\{x=0\}\backslash\{(0,0)\}$ and $\{y=0\}\backslash\{(0,0)\}$.
\end{example}

One of our main tools for proving results about the Minkowski sum is an alternative description of it in terms of the join of the two varieties.

For $X,Y$ subvarieties of $\mathbb{A}^n$ or $\mathbb{P}^n$, we let $J_{\mathrm{set}}(X,Y)$ be the {\it setwise join} of $X$ and $Y$, i.e., the union of the lines connecting distinct points $x\in X$ and $y\in Y$. This space is usually not closed and its Zariski closure $J(X,Y)$ is the classical {\it join} of $X$ and $Y$.

Our analysis of the Minkowski sum of affine algebraic sets $X$ and $Y$ via a join will involve hyperplanes positioned as in Lemma \ref{lemma:JoinHyperplaneIntersection} below. For an intuitive sense of the statement of the lemma, one may consider the case where $L,M,\text{ and }N$ are the projectivizations of parallel affine hyperplanes. %We will primarily focus on the case where projectivizing $X$ and $Y$ does not change their intersection.
%Insert/replace this last sentence with: Interested in when the setwise Minkowski sum is closed. Recall example of 2 hyperbolas. In fact, every situation in which the setwise Minkowski sum is not closed must look something like this.

\begin{lemma}\label{lemma:JoinHyperplaneIntersection}
Let $L,M,N$ be three distinct hyperplanes in $\PP^n$ with $E:=L\cap M=L\cap N=M\cap N$. Say $X\subset M$ and $Y\subset N$ are nonempty disjoint subvarieties. Let $X^a=X\setminus E$, $Y^a=Y\setminus E$, $\partial X=X\cap E$, and $\partial Y=Y\cap E$. Then:
\begin{enumerate}
\item[(i)] $J(X,Y)=J_{\mathrm{set}}(X,Y)$,
\item[(ii)] $J(X,Y)\cap L=(J_{\mathrm{set}}(X^a,Y^a)\cap L)\cup J_{\mathrm{set}}(\partial X,\partial Y)\cup\partial X\cup \partial Y$, and
\item[(iii)] $J(X,Y)\cap L\setminus E=J_{\mathrm{set}}(X^a,Y^a)\cap L$.
\end{enumerate}
 In particular, if $X$ and $Y$ have positive dimension then $$J(X,Y)\cap L=(J_{\mathrm{set}}(X^a,Y^a)\cap L)\cup J_{\mathrm{set}}(\partial X,\partial Y).$$
\end{lemma}\begin{proof}
\noindent {\it (i)} Because $X$ and $Y$ are disjoint, we have $J_{\mathrm{set}}(X,Y)$ is Zariski closed, so $J(X,Y)=J_{\mathrm{set}}(X,Y)$ [see Example 6.17 on p.70 of \cite{Harris}].

\smallskip
\noindent {\it (ii)} From the first part, we have 
$$J(X,Y)=J_{\mathrm{set}}(X^a,Y^a)\cup J_{\mathrm{set}}(X^a,\partial Y)\cup J_{\mathrm{set}}(\partial X,Y^a)\cup J_{\mathrm{set}}(\partial X,\partial Y).$$
So to get the claimed expression for $J(X,Y)\cap L$ it suffices to show 
\begin{enumerate}
\item[(a)]$J_{\mathrm{set}}(X^a,\partial Y)\cap L,J_{\mathrm{set}}(\partial X,Y^a)\cap L\subset \partial X\cup\partial Y$ and
\item[(b)] $J_{\mathrm{set}}(\partial X,\partial Y)\cup\partial X\cup \partial Y\subset L$.
\end{enumerate}
\begin{enumerate}
\item[(a)] By symmetry it is enough to show that $J_{\mathrm{set}}(X^a,\partial Y)\cap L\subset \partial Y$. 

Say $x\in X^a$ and $y\in\partial Y$. So $y\in L$ but $x\notin L$. Thus, the line between $x$ and $y$ intersects $L$ in exactly $\{y\}\subset\partial Y$.

\item[(b)] We show that $J_{\mathrm{set}}(\partial X,\partial Y)\cup\partial X\cup \partial Y\subset E$. 

By definition, $\partial X,\partial Y\subset E$. So, because $E$ is a linear space, for any $x\in\partial X$ and $y\in\partial Y$, the line between $x$ and $y$ is contained in $E$.
\end{enumerate}

\noindent {\it (iii)} First, note that because $J_{\mathrm{set}}(\partial X,\partial Y)\cup\partial X\cup \partial Y\subset E$, we have 
$$\left(J(X,Y)\cap L\right)\setminus E\subset \left(J_{\mathrm{set}}(X^a,Y^a)\cap L\right)\setminus E.$$
Hence, we just need to show that $J_{\mathrm{set}}(X^a,Y^a)\cap L$ is disjoint from $E$. 

Considering any $x\in X^a$ and $y\in Y^a$, it suffices to show that the line $\ell$ between $x$ and $y$ does not meet $E$. If we assume, towards a contradiction, that there is some $z\in\ell\cap E$, then $z$ and $x$ would be distinct points on the hyperplane $M$, so the line $\ell$ between them would be contained in $M$. But $y \in Y^a\subset N\setminus E=N\setminus M$, so $\ell$ cannot be contained in $M$.

Finally, if $X$ and $Y$ are positive dimensional then $\partial X=X\cap L$ and $\partial Y=Y\cap L$ are nonempty, so $\partial X,\partial Y\subset J_{\mathrm{set}}(\partial X,\partial Y)$.\qed
\end{proof}

Our alternative description of the Minkowski sum will give us cases in which $X+_{\mathrm{set}}Y$ is already closed. Recall from Example \ref{example:TwoHyperbolas} that for the two plane curves $X=\{xy=1\}$ and $Y=\{xy=-1\}$, $X+_{\mathrm{set}}Y$ is not Zariski closed. Note that in this example, $X$ and $Y$ have a common asymptote, or equivalently, that their projective closures meet at the line at infinity. We will see that when the characteristic of the base field is not 2, all cases where $X+_{\mathrm{set}}Y$ is not closed share an analogous property.

\begin{definition}
Let $X,Y \subset \mathbb{A}^n$ be varieties and denote the projective closures of $X$ and $Y$ in $\mathbb{P}^n$ by $\overline{X}$ and $\overline{Y}$, respectively. Let $H_0=\{[x_0:\cdots:x_n] \in \mathbb{P}^n \mid x_0=0\}$ be the {\it hyperplane at $\infty$}. %We say that $X \cap Y$ is \emph{finitely contained} if $\bar{X} \cap \bar{Y} \cap H = \emptyset$.
We say that $X$ and $Y$ are \emph{disjoint at infinity} if $\overline{X} \cap \overline{Y} \cap H_0 = \emptyset$.
\end{definition}

\begin{remark}
If $X$ and $Y$ are disjoint at infinity then $\dim(X\cap Y)<1$, thus $\dim X+\dim Y\leq n$.
\end{remark}

%\begin{proposition}\label{prop:CayleyTrick}
%Assume the characteristic of the base field is not 2.\\
%Let $X,Y\subset\mathbb{A}^n$ be algebraic sets. Let $\bar{X},\bar{Y}\subset\mathbb{P}^n$ denote the projective closures of $X,Y$, respectively, and let $H=\{[x_0:\cdots:x_n]\in\mathbb{P}^n|x_0=0\}$, ``the hyperplane at $\infty$''. Suppose that $\bar{X}\cap\bar{Y}\cap H=\emptyset$, i.e.\ the projective closures of $X$ and $Y$ don't intersect at $\infty$.\\
%Let $z_0,z_1$ be distinct scalars, let $\tilde{X},\tilde{Y}\subset\PP^{n+1}$ be the projective closures of $X\times\{z_0\}$ and $Y\times\{z_1\}$, respectively, and let $x_0,x_1,\ldots,x_n,z$ be the coordinates on $\PP^{n+1}$. We identify the hyperplane $S=\{z=\frac{z_0+z_1}{2}x_0\}\subset\PP^{n+1}$ with $\PP^{n}$, and so we identify $S\setminus H$ with $\mathbb{A}^n$. Then $J_{\mathrm{set}}(\tilde{X},\tilde{Y})=J(\tilde{X},\tilde{Y})$, $\frac{1}{2}(X+Y)=J(\tilde{X},\tilde{Y})\cap S\setminus H$, and $X+_{\mathrm{set}}Y=X+Y$, i.e.\ $X+_{\mathrm{set}}Y$ is Zariski closed.
%\end{proposition}
\begin{proposition}\label{prop:CayleyTrick}
Assume that the characteristic of $\Bbbk$ is not 2. 
Suppose $X, Y \subset\mathbb{A}^n$ are varieties that are disjoint at infinity. Let $z_0,z_1$ be distinct scalars and let $\tilde{X},\tilde{Y}\subset\PP^{n+1}$ be the projective closures of $X\times\{z_0\}$ and $Y\times\{z_1\}$, respectively. Let $x_0,x_1,\ldots,x_n,z$ be the coordinates on $\PP^{n+1}$. 

If we identify $S=\{z=\frac{z_0+z_1}{2}x_0\}\subset\PP^{n+1}$ with $\PP^{n}$ and $S\setminus H_0$ with $\mathbb{A}^n$, then: 
\begin{enumerate}
\item[(i)] $J_{\mathrm{set}}(\tilde{X},\tilde{Y})=J(\tilde{X},\tilde{Y})$;
\item[(ii)] $\frac{1}{2}(X+Y)=J(\tilde{X},\tilde{Y})\cap S\setminus H_0$;
\item[(iii)] $X+_{\mathrm{set}}Y=X+Y$, namely, $X+_{\mathrm{set}}Y$ is Zariski closed.
\end{enumerate}

\end{proposition}\begin{proof}
{\it (i)} Let $E=\{x_0=0\}\subset\mathbb{P}^{n+1}$ be the hyperplane at $\infty$ in $\mathbb{P}^{n+1}$. 
%Line 288 is rephrasing part of the original proof to use new term. Original is line 289
Note that $$E\cap\{z=z_0x_0\} = E\cap\left\{z=\frac{z_0+z_1}{2}x_0\right\} = E\cap\{z=z_1x_0\} = \{z=0,x_0=0\},$$ which is identified with $H_0$. Therefore the statement that $X$ and $Y$ are disjoint at infinity is equivalent to
$$\tilde{X} \cap \tilde{Y} \cap E=\emptyset.$$
%Note that because $\bar{X}\cap\bar{Y}\cap H=\emptyset$, (and because $E\cap\{z=z_0x_0\},E\cap\{z=\frac{z_0+z_1}{2}x_0\},$ and $E\cap\{z=z_1x_0\}$ all equal $\{z=0,x_0=0\}$ which is identified with $H$,) $\tilde{X}\cap\tilde{Y}\cap E=\emptyset$. 
On the other hand, $\tilde{X}\setminus E=X\times\{z_0\}$ and $\tilde{Y}\setminus E=Y\times\{z_1\}$, and so we see that $\tilde{X}\cap\tilde{Y}=\emptyset$. So, by Lemma \ref{lemma:JoinHyperplaneIntersection} applied to $S=\{z=\frac{z_0+z_1}{2}x_0\}$, $\tilde{X}\subset\{z=z_0x_0\}$, and $\tilde{Y}\subset\{z=z_1x_0\}$, we find that 
$$J_{\mathrm{set}}(\tilde{X},\tilde{Y})=J(\tilde{X},\tilde{Y})~~\text{ and }~~J(\tilde{X},\tilde{Y})\cap S\setminus H_0=J_{\mathrm{set}}(X\times\{z_0\},Y\times\{z_1\})\cap S.$$
\noindent{\it (ii)\ \& (iii)} For any $x\in X$ and $y\in Y$, the line between the points $(x,z_0)\in X\times\{z_0\}$ and $(y,z_1)\in Y\times\{z_1\}$ meets the affine hyperplane $S\setminus E=\{z=\frac{z_0+z_1}{2}\}$ in exactly the point $(\frac{x+y}{2},\frac{z_0+z_1}{2})$. So, we have shown that $J(\tilde{X},\tilde{Y})\cap S\setminus H_0=\frac{1}{2}(X+_{\mathrm{set}}Y)$. In particular, because $J(\tilde{X},\tilde{Y})$ is closed, this tells us that $X+_{\mathrm{set}}Y$ is a closed subset of $S\setminus H_0\cong\mathbb{A}^n$. Hence, $\frac{1}{2}(X+Y)=\frac{1}{2}(X+_{\mathrm{set}}Y)=J(\tilde{X},\tilde{Y})\cap S\setminus H_0$.
%
%\smallskip
%\textcolor{blue}{{\it (iii)} It follows from the proof of point $(ii)$.}
\qed

\end{proof}

\begin{remark}\label{remark:Cayley}We call the construction $\frac{1}{2}(X+Y)=J(\tilde{X},\tilde{Y})\cap S\setminus H_0$ the {\it Cayley trick}, as the underlying idea is exactly the same as that of the Cayley trick used to construct Minkowski sums of polytopes.
\end{remark}

 As a consequence of the following lemma, if we restrict to the cases with $\dim X+\dim Y\leq n$, then the hypothesis that $X$ and $Y$ are disjoint at infinity is a {\it genericity condition}. 

\begin{lemma}\label{lemma:GenericTransformationIntersection}
Let $X,Y\subset\PP^n$ be varieties with $\dim X+\dim Y<n$. Then, we have that the set $\{g\in GL_{n+1}~|~gX\cap Y=\emptyset\}$ is a nonempty open subset of $GL_{n+1}$. That is, for generic $g\in GL_{n+1}$, $gX$ and $Y$ do not intersect. 
\end{lemma}\begin{proof}
First, note that for any point $p\in\PP^n$ the stabilizer of $p$ in $GL_{n+1}$ has dimension $n^2+n+1$. This is because any two point stabilizers in $GL_{n+1}$ are conjugate and the stabilizer of $[1:0:0:\cdots:0]\in\PP^n$ is the set of all $g\in GL_{n+1}$ with first column of the form $\begin{bmatrix}*&0&0&\cdots&0\end{bmatrix}^T$, which has dimension $(n+1)n+1$.

Let $Z=\{(g,x,y)\in GL_{n+1}\times X\times Y~|~gx=y\}$ which is a subvariety of $GL_{n+1}\times X\times Y$. Let $\pi_1:Z\to GL_{n+1}$ and $\pi_2:Z\to X\times Y$ be the restrictions of the canonical projections from $GL_{n+1}\times X\times Y$. Note that $\pi_2$ is surjective, because for any $x,y\in\PP^n$ there exists some $g\in GL_{n+1}$ taking $x$ to $y$. Further, the fiber over any point $(x,y)\in X\times Y$ is a left coset of a point stabilizer in $GL_{n+1}$ and so has dimension $n^2+n+1$. Thus, $\dim Z=n^2+n+1+\dim X+\dim Y<n^2+2n+1=\dim GL_{n+1}$.

Because $X\times Y$ is projective, the projection $GL_{n+1}\times X\times Y\to GL_{n+1}$ is a closed map, so $\pi_1(Z)$ is a closed subset of $GL_{n+1}$. Since $\dim \pi_1(Z)\leq\dim Z<\dim GL_{n+1}$, $\pi_1(Z)$ is a proper closed subset of $GL_{n+1}$. So, $$\{g\in GL_{n+1}|gX\cap Y=\emptyset\}=GL_{n+1}\setminus\pi_1(Z)$$ is a nonempty open subset of $GL_{n+1}$.
\qed
\end{proof}

Now, we claim that if $X, Y \subset\A^n$ are varieties with $\dim X + \dim Y \leq n$, then
\begin{center}
{\it for general $g\in GL_{n}$, $gX$ and $Y$ are disjoint at infinity.}
\end{center}
To see this, note that, considering $H_0 = \{x_0 = 0\}$, $$\dim (\overline{X}\cap H_0)+\dim(\overline{Y}\cap H_0)\leq\dim X-1+\dim Y-1<\dim X+\dim Y-1\leq n-1.$$
%$\Big($So, if we identify $H_0$ with \textcolor{green}{$\PP(\A^n)$}\marginpar{why are we using this symbol here?} so that $GL_n$ acts on $H_0$, then for $g\in GL_n$, we have $\Big)$
The action of $GL_n$ on $\A^n$ extends to an action on $\PP^n$, and the identification $H_0\cong\PP^{n-1}$ is $GL_n$-equivariant. So we have

$$\overline{gX}\cap\overline{Y}\cap H_0=(\overline{gX}\cap H_0)\cap(\overline{Y}\cap H_0)=g(\overline{X}\cap H_0)\cap(\overline{Y}\cap H_0).$$ By Lemma \ref{lemma:GenericTransformationIntersection}, for general $g\in GL_n$ this is empty. 

\begin{remark}
We could have used the group of affine transformations $\mathrm{Aff}_n=\A^n\rtimes GL_n$. Indeed, shifting an affine variety does not change the part at infinity of its projective closure.
\end{remark}

When a result holds under the same conditions as Lemma \ref{lemma:GenericTransformationIntersection}, i.e., if we fix $X$ and $Y$ then it holds for $gX$ and $Y$, for general $g\in\mathrm{Aff}_n$, we shall say that the result holds for $X$ and $Y$ \textit{in general position}.%``in general position''.}

%As a kind of converse to this, we have that if $X$ and $Y$ intersect boundedly then $\dim X+\dim Y\leq n$. For if $\dim X+\dim Y>n$ then $\dim(\bar{X}\cap \bar{Y})\geq\dim X+\dim Y-n>0$ so $\bar{X}\cap\bar{Y}\cap H_0\neq\emptyset$.

%We are now ready to compute the dimension of Minkowski sums. Based on small examples as in Section \ref{sec:experiments}, it may seem that for $\dim X+\dim Y\leq n$, except for some very special cases we have $\dim(X+Y)=\dim X+\dim Y$. As we prove in the following, the condition of {\it bounded intersection} is sufficient to have additivity. 
We are now ready to compute the dimension of Minkowski sums. Based on the examples in Section \ref{sec:experiments}, it seems that for $\dim X+\dim Y\leq n$, we get $\dim(X+Y)=\dim X+\dim Y$. This does happen generically.
%As we prove in the following, the condition of {\it bounded intersection} is sufficient to have additivity. 

\begin{theorem}\label{thm:CharIndepMinkowskiDimension}
%Let $X,Y\subset\A^n$ be varieties with $\dim X+\dim Y\leq n$. Then for $X$ and $Y$ in general position, $\dim(X+Y)=\dim X+\dim Y$.
Let $X,Y\subset\A^n$ be varieties. Then for $X$ and $Y$ in general position, $\dim(X+Y)=\min\{\dim X+\dim Y,n\}$.
\end{theorem}\begin{proof}
As observed in Note \ref{note:SumAndProductBasicProperties}, we have that, for any $X,Y \subset \A^n$,
$$\dim(X+Y)\leq\min\{\dim(X)+\dim(Y),n\}.$$
So, we need to prove the converse for $X,Y$ in general position. 

Let $k=\dim(X)$ and $l=\dim(Y)$.

Note that because $(X+v)+Y=(X+Y)+v$ for any vector $v$, it suffices to show that for general $g\in GL_n$, $\dim(gX+Y)\geq\min\{\dim(X)+\dim(Y),n\}$. We consider the case $\dim(X)+\dim(Y)\leq n$, the case $\dim(X)+\dim(Y)\geq n$ being analogous.

Note that by just looking at full-dimensional irreducible components of $X$ and $Y$, we may assume without loss of generality that $X$ and $Y$ are irreducible. 

We denote by $T_p(X)$ the {\it tangent space} to the variety $X$ at the point $p$.

For now fix $g\in GL_n$. If $(p,q)\in\A^n\times\A^n$ then 
$$(d\phi_+)_{(p,q)}:T_p\A^n\times T_q\A^n\to T_{p+q}\A^n$$
 is simply the addition map $\phi_+$, and so we see that if $p\in gX$ and $q\in Y$ then $T_p(gX)+T_qY\subseteq T_{p+q}(gX+Y)$. So to conclude that $\dim(gX+Y)\geq\dim(X)+\dim(Y)$ it suffices to show that there is a dense subset $\Xi$ of $gX+Y$ such that for each $\xi\in\Xi$ there exist $p\in gX$ and $q\in Y$ with $\xi=p+q$ and $T_p(gX)\cap T_qY=0$, for then 
\begin{align*}
\dim(T_\xi(gX+Y))&\geq\dim(T_p(gX)+T_qY)\\
&=\dim(T_p(gX))+\dim(T_qY)\geq\dim(X)+\dim(Y),
\end{align*}
 and because $\Xi$ is dense some $\xi\in\Xi$ is a smooth point of $gX+Y$. Also, because the image of a dense subset under a continuous function is a dense subset of the image, we see that it suffices to show that there is a nonempty open subset of $gX\times Y$ such that for $(p,q)$ in this set, $T_p(gX)\cap T_qY=0$.

For any variety $Z\subset\A^n$ let $Z_{\mathrm{sm}}$ denote the smooth locus of $Z$. So we have the morphism $\psi_Z:Z_{\mathrm{sm}}\to Gr(\dim(Z),\A^n)$, $p\mapsto T_pZ$ and we let $\Psi_Z$ denote the image of this morphism.

Consider 
$$U=\{(V,W)\in Gr(k,\A^n)\times Gr(l,\A^n)|V\cap W\neq0\},$$
 which is an open subset of $Gr(k,\A^n)\times Gr(l,\A^n)$. In particular, if we let 
$$\varphi_g:=\psi_{gX}\times\psi_Y:gX_{\mathrm{sm}}\times Y_{\mathrm{sm}}\to Gr(k,\A^n)\times Gr(l,\A^n),$$
 then $\varphi_g^{-1}(U)$ is a (possibly empty) open subset of $gX\times Y$, and if $(p,q)\in\varphi_g^{-1}(U)$ then $T_p(gX)\cap T_qY=0$. Also, $\varphi_g^{-1}(U)$ is nonempty if and only if $(\Psi_{gX}\times\Psi_{Y})\cap U$ is nonempty. So we conclude that to show that $\dim(gX+Y)\geq\dim(X)+\dim(Y)$, it suffices to show that $(\Psi_{gX}\times\Psi_Y)\cap U\neq\emptyset$.

Now we let $g\in GL_n$ vary. Fix $p\in X_{\mathrm{sm}}$ and $q\in Y_{\mathrm{sm}}$. So for $g\in GL_n$, $gp\in (gX)_\mathrm{sm}$ with $T_{gp}(gX)=g(T_pX)$. Now $T_qY$ is an $l$-dimensional subspace of $\A^n$ and so because $k+l\leq n$, $\{V\in Gr(k,\A^n)|V\cap T_qY=0\}$ is a nonempty open subset of the Grassmannian. So because $GL_n$ acts transitively on $Gr(k,\A^n)$ we conclude that for generic $g\in GL_n$, $T_{gp}(gX)\cap T_qY=g(T_pX)\cap T_qY=0$. Thus $(gp,q)\in(\Psi_{gX}\times\Psi_Y)\cap U$, and so $\dim(gX+Y)\geq\dim(X)+\dim(Y)$.

For the case where $\dim(X)+\dim(Y)\geq n$ the same proof works upon replacing the condition that tangent spaces intersect trivially with the condition that they intersect transversely.
\qed
\end{proof}

Further, when the characteristic of the base field is not 2 we can use the Cayley trick to show that the condition of {\it disjoint at infinity} is sufficient to have additivity of dimension.

\begin{theorem}\label{thm:MinkowskiDimension}
Assume the characteristic of the base field is not 2. 
Let $X,Y\subset\A^n$ be varieties which are disjoint at infinity. Then $\dim(X+Y)=\dim X+\dim Y$.
%If the projective closures of $X$ and $Y$ don't meet at infinity then $\dim(X+Y)=\dim X+\dim Y$.
\end{theorem}\begin{proof}
As observed in Note \ref{note:SumAndProductBasicProperties}, we have that for any $X,Y \subset \A^n$,
$$\dim(X+Y)\leq\dim(X)+\dim(Y).$$ If either $X$ or $Y$ has dimension zero then $X+Y$ is a union of finitely many shifts of the other and so has dimension $\dim X+\dim Y$. 

Assume $X$ and $Y$ have both positive dimension. Then, by Proposition \ref{prop:CayleyTrick} (with any $z_0,z_1$) and Lemma \ref{lemma:JoinHyperplaneIntersection}, we have that $$J(\tilde{X},\tilde{Y})\cap S=\frac{1}{2}(X+Y)\cup J(\partial X,\partial Y)$$ where $\partial X=\overline{X}\cap H_0$ and $\partial Y=\overline{Y}\cap H_0$ are the parts at infinity of the projective closures of $X$ and $Y$, and $\frac{1}{2}(X+Y)$ is an open subset of $J(\tilde{X},\tilde{Y})$ while $J(\partial X,\partial Y)$ is closed. Hence, we have
\begin{align*}
\dim X+\dim Y&=\dim J(\tilde{X},\tilde{Y})-1\leq\dim(J(\tilde{X},\tilde{Y})\cap S)\\
&=\dim \left(\frac{1}{2}(X+Y)\cup J(\partial X,\partial Y)\right)\\
&=\max\left\{\dim(X+Y),\dim J(\partial X,\partial Y)\right\}\\
&=\max\left\{\dim(X+Y),\dim X+\dim Y-1\right\},
\end{align*}
where last equality follows since$$\dim J(\partial X,\partial Y)=\dim X-1+\dim Y-1+1=\dim X+\dim Y-1.$$
So $\dim(X+Y)=\dim X+\dim Y$.
\qed
\end{proof}

%{\color{red}(Insert discussion about: Seemed like maybe being badly ruled would be the right condition, but this isn't right. Example of $I_1=(x_1^2-x_2,x_4^3+x_5^3+x_6^3,x_3)$ and $I_2=(y_2^2-y_3,y_4^3+y_5^2+y_6^3,y_1)$ Having Minkowski sum with ideal $(z_1^4 - 2*z_1^2*z_2 + z_2^2 - z_3)$. Explain idea where this example comes from.)}

We now consider the degree of Minkowski sums. Recall that the \textit{degree} of a variety $X$ of dimension $d$ in $\A^n$ or $\PP^n$ is the number of points in the intersection of $X$ and a general linear subspace of dimension $n-d$.

\begin{proposition}\label{prop:DegreeForMinkowskiSumWithBoundedIntersection}
Let $\Bbbk$ be the ground field with characteristic other than 2. 
% Let $X,Y\subset\A^n$ be algebraic sets with $\dim X+\dim Y\leq n$, and suppose that $X$ and $Y$ intersect boundedly.
Let $X,Y\subset\A^n$ be varieties which are disjoint at infinity.
Then, for generic $\alpha\in \Bbbk^\times$, in the same notation as in Proposition \ref{prop:CayleyTrick}, we have that $\deg(\alpha X+Y)=\deg \left(J(\tilde{X},\tilde{Y})\right)$.
\end{proposition}\begin{proof}
The proof will go in three main steps.

\begin{enumerate}
\item[(i)] Show that, up to projective equivalence, dilating $X$ by a generic $\alpha\in \Bbbk^\times$ and then applying the Cayley trick is the same as intersecting $J(\tilde{X},\tilde{Y})$ with a generic hyperplane whose affine part is parallel to $S\backslash H_0$.
\item[(ii)] Prove that for generic $\alpha$ the corresponding hyperplane intersects $J(\tilde{X},\tilde{Y})$ generically transversely.
\item[(iii)] Apply B\'{e}zout's theorem and show that the part of the intersection that is at infinity does not contribute to the degree.
\end{enumerate}
Once again we use our Cayley trick and, to simplify computations, we fix $z_0=0$ and $z_1=1$. Note that, for any $\alpha\in \Bbbk^\times$, $\alpha X$ and $Y$ are disjoint at infinity
%the projective closures of $\alpha X$ and $Y$ still don't meet at infinity, 
so we still get the conclusions of Proposition \ref{prop:CayleyTrick} and Theorem \ref{thm:MinkowskiDimension}.

\smallskip
\noindent (i) For $\alpha\in \Bbbk^\times$, let $$\Phi_\alpha=\begin{pmatrix}
1&0 ~~\cdots~~ 0&\alpha-1\\
\begin{matrix} 0 \\ \vdots \\ 0 \end{matrix}&\alpha I_n&\begin{matrix} 0 \\ \vdots \\ 0 \end{matrix}\\
0& 0 ~~ \cdots ~~ 0 &\alpha
\end{pmatrix}\in GL_{n+2}.$$ We consider $GL_{n+2}$ as acting on $\PP^{n+1}$ with coordinates $x_0,x_1,\ldots,x_n,z$. 
For $\alpha,\beta\in \Bbbk^\times$ we have $$\Phi_\alpha \Phi_\beta=\begin{pmatrix}
1&&\alpha-1\\
&\alpha I_n\\
&&\alpha
\end{pmatrix}\begin{pmatrix}
1&&\beta-1\\
&\beta I_n\\
&&\beta
\end{pmatrix}=\begin{pmatrix}
1&&\beta-1+\beta(\alpha-1)\\
&\alpha\beta I_n\\
&&\alpha\beta
\end{pmatrix}=\Phi_{\alpha\beta},$$so $\alpha\mapsto \Phi_\alpha$ is a group homomorphism $\Bbbk^\times\to GL_{n+2}$.

Note that $\Phi_\alpha$ acts on the hyperplane $\{z=0\}$ as 
$$\Phi_\alpha\begin{bmatrix}1\\\underline{x}\\0\end{bmatrix}
=\begin{bmatrix}1\\\alpha\underline{x}\\0\end{bmatrix}~~\text{ and }~~\Phi_\alpha\begin{bmatrix}0\\\underline{x}\\0\end{bmatrix}
=\begin{bmatrix}0\\\alpha\underline{x}\\0\end{bmatrix}
=\begin{bmatrix}0\\\underline{x}\\0\end{bmatrix}.$$ 
Similarly, $\Phi_\alpha$ fixes the hyperplane $\{z=x_0\}$ pointwise. Thus, $\Phi_\alpha(\tilde{X})=\widetilde{\alpha X}$ and $\Phi_\alpha(\tilde{Y})=\tilde{Y}$. Since $\Phi_\alpha$ acts as a projective transformation, and so takes lines to lines, $\Phi_\alpha(J(\tilde{X},\tilde{Y}))=J(\widetilde{\alpha X},\tilde{Y})$. In particular, $$\deg \left(J\left(\widetilde{\alpha X},\tilde{Y}\right)\right)=\deg \left(J(\tilde{X},\tilde{Y})\right).$$
We know that $\frac{1}{2}(\alpha X+Y)=J\left(\widetilde{\alpha X},\tilde{Y}\right)\cap S\setminus H_0$, so we consider $\Phi_\alpha^{-1}(S\setminus H_0)$. We get $$S\setminus H_0=\big\{z=\frac{1}{2}x_0\big\}\setminus\{z=0,x_0=0\}=\big\{x_0=1,z=\frac{1}{2}\big\},$$ and we find that, for $\alpha\neq-1$,
\begin{align*}
\Phi_\alpha^{-1}\begin{bmatrix} 1\\\underline{w}\\1/2\end{bmatrix}
&=\begin{pmatrix}
1&&\alpha^{-1}-1\\
&\alpha^{-1} I_n\\
&&\alpha^{-1}
\end{pmatrix}\begin{bmatrix} 1\\\underline{w}\\1/2\end{bmatrix}
=\begin{bmatrix}\frac{\alpha^{-1}+1}{2}\\\alpha^{-1}\underline{w}\\\alpha^{-1}/2\end{bmatrix}=\begin{bmatrix}1\\\frac{2\alpha^{-1}}{\alpha^{-1}+1}\underline{w}\\\frac{\alpha^{-1}}{\alpha^{-1}+1}\end{bmatrix}
=\begin{bmatrix}1\\\frac{2}{1+\alpha}\underline{w}\\\frac{1}{1+\alpha}\end{bmatrix}.
\end{align*}
Thus $\Phi_\alpha^{-1}(S\setminus H_0)=\{x_0=1,z=\frac{1}{1+\alpha}\}=\{z=\frac{1}{1+\alpha}x_0\}\setminus H_0$.
\smallskip

\noindent (ii) We claim that, for generic $\alpha$, 
\begin{center}
{\it $\left\{z=\frac{1}{1+\alpha}x_0\right\}$ intersects $J(\tilde{X},\tilde{Y})$ generically transversely.}
\end{center}
 First, note that it suffices to only consider the {\it affine points} of $J(\tilde{X},\tilde{Y})$, i.e.\ those with $x_0=1$, because
 \begin{align*}\dim\left(J(\tilde{X},\tilde{Y})\cap \left\{z=\frac{1}{1+\alpha}x_0\right\}\setminus H_0\right) & =\dim\left(J\left(\widetilde{\alpha X},\widetilde{Y}\right)\cap S\setminus H_0\right) \\ &=\dim(\alpha X+Y)=\dim X+\dim Y \\ &> \dim J(\partial X,\partial Y) =\dim(J(\tilde{X},\tilde{Y})\cap H_0),
 \end{align*}
  But $\A^{n+1}\subset\PP^{n+1}$ is the disjoint union of $\{x_0=1,z=a\}$ as $a$ ranges over $\Bbbk$, so for all but finitely many $a\in k$, $\{x_0=1,z=a\}\cap J(\tilde{X},\tilde{Y})$ must not be contained in the singular locus of $J(\tilde{X},\tilde{Y})$. So, for all but finitely many $\alpha\in \Bbbk^\times$, we have that $\{x_0=1,z=\frac{1}{1+\alpha}\}\cap J(\tilde{X},\tilde{Y})$ must not be contained in the singular locus of $J(\tilde{X},\tilde{Y})$. So for generic $\alpha\in \Bbbk^\times$, the general point of $\{z=\frac{1}{1+\alpha}x_0\}\cap J(\tilde{X},\tilde{Y})$ is a smooth point of $J(\tilde{X},\tilde{Y})$. In order to check transversality, we need another description of this intersection, which we compute now.
\begin{align*}
J(\tilde{X},\tilde{Y})\cap\left\{z=\frac{1}{1+\alpha}\right\}&=\Phi_\alpha^{-1}\left(J\left(\widetilde{\alpha X},\tilde{Y}\right)\cap S\setminus H_0\right)\\
&=\Phi_\alpha^{-1}\left(J(\alpha X\times\{0\},Y\times\{1\})\cap S\setminus H_0\right)\\
&=J\left(X\times\{0\},Y\times\{1\}\right)\cap \left\{z=\frac{1}{1+\alpha}\right\}
\end{align*}
 where the second equality follows from Lemma \ref{lemma:JoinHyperplaneIntersection}. 
 
 Thus, considering $p\in J(\tilde{X},\tilde{Y})\cap \{z=\frac{1}{1+\alpha}\}$, we have that $p$ is on the line between the points $(x,0)$ and $(y,1)$, for some $x\in X$ and $y\in Y$. Since this line intersects $S=\{z=\frac{1}{1+\alpha}x_0\}$ transversely and $T_pJ(\tilde{X},\tilde{Y})$ contains this line, if $p$ is a smooth point of $J(\tilde{X},\tilde{Y})$ then we have that $\{z=\frac{1}{1+\alpha}x_0\}$ and $J(\tilde{X},\tilde{Y})$ intersect transversely at $p$. Thus, for generic $\alpha$, $\{z=\frac{1}{1+\alpha}x_0\}$ intersects $J(\tilde{X},\tilde{Y})$ generically transversely.

\smallskip
\noindent (iii) For such an $\alpha$, applying B\'{e}zout's theorem gives us that 
$$\deg\left(J\left(\tilde{X},\tilde{Y}\right)\cap \left\{z=\frac{1}{1+\alpha}x_0\right\}\right)=\deg \left(J(\tilde{X},\tilde{Y})\right).$$
We can write $J(\tilde{X},\tilde{Y})\cap \{z=\frac{1}{1+\alpha}x_0\}$ as the disjoint union of the open subset $J(\tilde{X},\tilde{Y})\cap\{x_0=1,z=\frac{1}{1+\alpha}\}$ and the closed subset $J(\tilde{X},\tilde{Y})\cap H_0$. Now,
\begin{align*}
J(\tilde{X},\tilde{Y})\cap\left\{x_0=1,z=\frac{1}{1+\alpha}\right\} & =\Phi_\alpha^{-1}\left(J\left(\widetilde{\alpha X},\tilde{Y}\right)\cap S\setminus H_0\right) \\ & =\Phi_\alpha^{-1}\left(\frac{1}{2}(\alpha X+Y)\right)
\end{align*}
 has dimension $\dim X+\dim Y$ and 
$J(\tilde{X},\tilde{Y})\cap H_0=J(\partial X,\partial Y)$
 has dimension $\dim X+\dim Y-1$. Therefore,
\begin{align*}\deg\left(J\left(\tilde{X},\tilde{Y}\right)\cap\left\{x_0=1,z=\frac{1}{1+\alpha}\right\}\right)&=\deg\left(J(\tilde{X},\tilde{Y})\cap \left\{z=\frac{1}{1+\alpha}x_0\right\}\right)\\ & = \deg \left(J(\tilde{X},\tilde{Y})\right).
\end{align*}
Finally, since 
$$\frac{1}{2}(\alpha X+Y)=\Phi_\alpha\left(J\left(\tilde{X},\tilde{Y}\right)\cap\left\{x_0=1,z=\frac{1}{1+\alpha}\right\}\right),$$
 we get $\deg(\alpha X+Y)=\deg \left(J(\tilde{X},\tilde{Y})\right)$.
\qed
\end{proof}

\begin{corollary}\label{coro:ActualDegreeComputation}
Suppose $\Bbbk$ has characteristic other than $2$. Let $X,Y\subset \A^n$ be varieties whose projective closures $\overline{X},\overline{Y}\subset\PP^n$ are contained in complementary linear subspaces; equivalently, $X,Y$ are contained in disjoint affine subspaces which are not parallel. Then for generic $\alpha\in\Bbbk^\times$, $\deg(\alpha X+Y) = \deg(X)\deg(Y)$.
\end{corollary}\begin{proof}
Since $\overline{X}$ and $\overline{Y}$ are contained in complementary linear spaces they are disjoint, so, in particular, $X$ and $Y$ are disjoint at infinity. 

By Proposition \ref{prop:DegreeForMinkowskiSumWithBoundedIntersection}, for a generic $\alpha\in \Bbbk^\times$, we have $\deg(\alpha X+Y)=\deg \left(J(\tilde{X},\tilde{Y})\right)$. Moreover, $\overline{X},\overline{Y}$ contained in complementary linear spaces also gives us that $\tilde{X}$ and $\tilde{Y}$ are contained in complementary linear spaces, so $\deg \left(J(\tilde{X},\tilde{Y})\right)=\deg(\tilde{X})\deg(\tilde{Y})$; see \cite[Example 18.17]{Harris}. So, for generic $\alpha\in \Bbbk^\times$, 
\begin{center}
$\deg(\alpha X+Y)=\deg \left(J(\tilde{X},\tilde{Y})\right)=\deg(\tilde{X})\deg(\tilde{Y})=\deg(X)\deg(Y).$
\end{center} \qed
\end{proof}

%{\color{orange} WHOOPS! Having just written out this whole proof, I just now realize that it is also true for a more trivial reason: If $\bar{X},\bar{Y}$ are in complementary linear spaces, then in the right basis, $X+Y$ is exactly the product of $X$ and $Y$ (I'm pretty sure). Should we say something about this? Should we take out this corollary?}

%\subsection*{Hadamard products of projective varieties}
\section{Hadamard products of projective varieties}\label{sec:Hadamard}

%{\color{red} If all Hadamard results are spoken of in terms of projective varieties, maybe we should reduce Hadamard product of of affine varieties to a footnote (e.g. "it is defined similarly for affine varieties and it is and the affine version of the following results hold by affine versions of the same arguments")?}

%\subsection{General facts}
We defined the Hadamard product of projective varieties $X,Y \subset \PP^n$ as 
$$
 X\star Y := \overline{\{p\star q ~:~ p\in X,~q\in Y,~p\star q \text{ is defined}\}}\subset \PP^n,
$$
where $p\star q$ is the point obtained by entry-wise multiplication of the points $p,q$. 

Also in this case the operation of closure is crucial.
\begin{example}\label{example:Hadamard_closure}
Consider the Hadamard product between the rational normal curve $\mathcal{C}_3 = \{[a^3:a^2b:ab^2:b^3] ~|~ [a:b]\in\PP^1\}$ in $\PP^3$ and the point $P = [0:1:1:0]$. Now, we obviously have $\mathcal{C}_3 \star P \subset \{z_0 = z_3 = 0\}$. The equality follows because, if $ab \neq 0$, then we have that $[0:a:b:0] = [a^3:a^2b:ab^2:b^3] \star [0:1:1:0]$. However, in this case the operation of taking the closure is needed in order to get the entire line; indeed, the points $[0:1:0:0]$ and $[0:0:1:0]$ cannot be written as the Hadamard product of a point in $\mathcal{C}_3$ and the point $P$.
\end{example}
Another useful way to describe the Hadamard product of projective varieties is as a linear projection of the Segre product of $X$ and $Y$, i.e., the variety obtained as the image of $X\times Y$ under the map
\begin{eqnarray*}
  \psi_{n,n}: & \PP^n \times \PP^n & \longrightarrow ~~~~~~~~~~~~~~~~~~~~~~~~~ \PP^{n^2+2n}, \\
  & ([a_0:\ldots:a_n],[b_0:\ldots:b_n]) &\mapsto [a_0b_0:a_0b_1:a_0b_2:\ldots:a_nb_{n-1}:a_nb_n].
\end{eqnarray*}
If $z_{ij}$, with $i = 0,\ldots,n$, $j = 0,\ldots,n$, are the coordinates of the ambient space of the Segre product $\PP^{n^2+2n}$, then the Hadamard product $X\star Y$ is the projection of $X\times Y$ with respect to the linear space $\{z_{ii} = 0 ~|~ i = 0,\ldots,n\}$.

Therefore, as observed in Note \ref{note:SumAndProductBasicProperties}, if $X$ and $Y$ are irreducible, then $X\star Y$ is irreducible and the dimension of their Hadamard product is at most the sum of the dimensions of the original varieties, i.e., $
 \dim(X\star Y) \leq \dim(X)+\dim(Y).
$
\begin{example}
It is easy to find examples where equality does not hold. Actually, the dimension of the Hadamard product of two varieties can be arbitrary small. E.g., consider two skew lines in $\PP^3$ as $H_{01} = H_0 \cap H_1=\{[0:0:a:b] ~|~ [a:b] \in \PP^1\}$ and $H_{23} = H_2\cap H_3 = \{[c:d:0:0] ~|~ [c:d] \in \PP^1\}$. Then $H_{01} \star H_{23}$ is empty.
\end{example}
A classic approach to compute the dimension of projective varieties is to look at their tangent space. From now, we consider $\CC$ as the ground field in order to avoid fuzzy behaviors caused by positive characteristics or non algebraically closed fields. Also, this is the case we want to consider in our applications.

In the case of joins, there is a result by A. Terracini \cite{T11} which describes the tangent space of the join at a generic point in terms of the tangent spaces of the original varieties. In \cite{BCK}, the authors proved a version of this result for Hadamard products of projective varieties.

\begin{lemma}\label{lemma:Terracini}{\rm\cite[Lemma 2.12]{BCK}} Let $p\in X$ and $q\in Y$ be generic points, then the tangent space to the Hadamard product $X\star Y$ at the point $p\star q$ is given by
$$
 T_{p\star q}(X\star Y) = \left\langle p\star T_qY, T_pX \star q\right\rangle.
$$
\end{lemma} 
Another powerful tool to study Hadamard products of projective varieties is tropical geometry. In particular, we have the following relation. Since we are not using tropical geometry elsewhere, here we assume the reader to be familiar with the concept of {\it tropicalization of a variety}. For the inexperienced reader, we suggest to read \cite{MS} for an introduction of the topic.

\begin{proposition}{\rm \cite[Proposition 5.5.11]{MS}}
 Given two irreducible varieties $X,Y\subset\PP^n$, the tropicalization of the Hadamard product of $X$ and $Y$ is the Minkowski sum of their tropicalizations, i.e.,
 \begin{center}$
  {\rm trop}(X\star Y) = {\rm trop}(X) + {\rm trop}(Y).
 $\end{center}
\end{proposition}
Applying this result, in \cite{BCK}, the authors gave an upper-bound for the dimension of the Hadamard product of two varieties.
\begin{proposition}\label{prop:expdim}{\rm \cite[Proposition 5.4]{BCK}}
 Let $X,Y\subset \PP^n$ be irreducible varieties. Let $H\subset (\CC^*)^{n+1}/\CC^*$ be the maximal subtorus acting on both $X$ and $Y$ and let $G\subset (\CC^*)^{n+1}/\CC^*$ be the smallest subtorus having a coset containing $X$ and a coset containing $Y$. Then 
 \begin{center}$
  \dim(X\star Y) \leq \min\{\dim(X) + \dim(Y) - \dim(H), \dim(G)\}.
 $\end{center}
\end{proposition}
We call this upper bound {\it expected dimension} and denote it ${\rm exp}.\dim(X\star Y)$. However, this is not always the correct dimension. In \cite{BCK}, the authors present an example of a Hadamard product of two projective varieties with dimension strictly smaller than the expected dimension. 

\bigskip

From the definition of the Hadamard product of two varieties, it makes sense also to analyze self Hadamard products of a projective variety. We call them {\it Hadamard powers} of a projective variety.
\begin{definition}
 We define the {\it $s$-th Hadamard power} of a projective variety $X$ as 
 $$
  X^{\star s} := X\star X^{\star (s-1)},\text{ for } s\geq 0,
 $$
 where $X^{\star 0} := [1:\ldots:1].$
\end{definition}
%In the recent paper \cite{CMS}, the authors use the Hadamard product of varieties to give a geometric interpretation of some particular stochastic model used in machine learning and data mining. In \cite{BCK}, the authors begun to study geometric properties of Hadamard product, and powers, of varieties, with a particular focus on linear spaces. 
In general, a projective variety is not contained in its Hadamard powers. However, if ${\sf 1}_n = [1:\ldots:1] \in \PP^n$ lies in the variety $X$, we get the following chain of non necessary strict inclusions
\begin{equation}\label{eq:chain}
% X \subset X^{\star 2} \subset \ldots\subset X^{\star s} \subset \ldots \subset \PP^n.
 X \subset X^{\star 2} \subset \dotsb\subset X^{\star s} \subset \dotsb \subset \PP^n.
\end{equation}
Therefore, it becomes very natural to check if the Hadamard powers of a projective variety $X$ eventually fill the ambient space.
In general, the answer is no.
\begin{proposition}\label{prop:binomial ideal}
Let $X$ be a toric variety in $\PP^n$. Then, $X = X^{\star 2}.$
\end{proposition}
\begin{proof}
 Since any toric variety contains the point $[1:\ldots:1]$, it follows that $X \subset X^{\star 2}$. The other inclusion follows by applying Proposition \ref{prop:expdim} to the case $X = Y = H $.
\qed
\end{proof}
\begin{remark}
 Recently, C. Bocci and E. Carlini gave a necessary and sufficient condition for a plane irreducible curve $C \subset \PP^2$ to have its $t$-th Hadamard power equal to the curve itself. This result has been shared with us in private communication and will appear in \cite{BC}.
\end{remark}
\begin{remark}
Proposition \ref{prop:binomial ideal} can be proved directly by recalling that the ideals defining toric varieties are given by {\it binomial ideals}, namely ideals whose generators are differences of monomials as $f_{\alpha,\beta} = x^\alpha - x^\beta$, where $\alpha,\beta\in \NN^{n+1}$ and we use the multi-index notation $x^\alpha := x_0^{\alpha_0}\cdots x_n^{\alpha_n}$.
 
 Now, consider two points of $X$, $p = [p_0:\ldots:p_n]$ and $q = [q_0:\ldots:q_n]$. For any generator $f_{\alpha,\beta}$ of the ideal defining $X$, we have $p^\alpha - p^\beta = q^\alpha - q^\beta = 0$. Therefore, 
 \begin{align*}
  (p\star q)^\alpha - (p\star q)^\beta & = p^\alpha q^\alpha - p^\beta  q^\beta = p^\alpha q^\alpha - p^\alpha q^\beta + p^\alpha  q^\beta -p^\beta  q^\beta = \\
  & = p^\alpha  (q^\alpha -q^\beta) - q^\beta (p^\alpha - p^\beta) = 0;
 \end{align*}
 hence, $p\star q \in X$.
\end{remark}

\begin{remark}
 Given a projective variety $X\subset \PP^n$, the $s${\it -th secant variety} $\sigma_s(X)$ is the Zariski closure of the union of linear spaces spanned by $s$ points lying on $X$. This is a very classical object that has been studied since the second half of $19$-th century. In particular, we have a chain of non necessary strict inclusions given by 
 $$
  X\subset \sigma_2(X) \subset \ldots \subset \sigma_s(X) \subset\ldots\subset\PP^n.
 $$
 Therefore, we can ask if the secant varieties of a variety $X$ eventually fill the ambient space. It is not difficult to prove that the answer is no. Indeed, if $H$ is a linear space, then $\sigma_2(H) = H$ and, therefore, if $X$ is degenerate, i.e., it is contained in a proper linear subspace of $\PP^n$, then its secant varieties do not fill the ambient space.

 Hadamard powers of projective varieties may be viewed as the multiplicative version of the classical notion of secant varieties where instead of looking at the linear span of points lying on a variety we consider their Hadamard product. Moreover, by Proposition \ref{prop:binomial ideal}, we have that the role played by linear spaces in the case of secant varieties is taken by toric varieties in the case of Hadamard products.
\end{remark}

\begin{example}\label{example:rank1}
 A concrete example satisfying the assumptions of Proposition \ref{prop:binomial ideal} is the variety $X_1 \subset \PP({\rm Mat}_{m,n})$ of rank $1$ matrices of size $m\times n$. Indeed, it is generated by the $2 \times 2$ minors of the generic matrix $(z_{ij})^{j = 1,\ldots,n}_{i = 1,\ldots,m}$. Therefore, $X_1^{\star 2} = X_1$. This gives another proof of the well-known fact that the Hadamard product of two rank $1$ matrices is still of rank $1$.
\end{example}
 The latter example rises a very interesting question. 
 \begin{question}What if we consider matrices of rank higher than $1$? Can we decompose {\it all} matrices as Hadamard products of rank $r>1$ matrices?\end{question} The answer is positive, as we show in the following proposition.
 \begin{proposition}\label{lemma:max rank}
  Let $M$ be a matrix of size $m\times n$ and fix $2\leq r\leq \min\{m,n\}$. Then, $M$ can be written as the Hadamard product of at most $\left\lceil \frac{\min\{m,n\}}{r-1}\right\rceil$ matrices of rank less or equal than $r$.
 \end{proposition}
 \begin{proof}
  Without loss of generality, we may assume that $m \leq n$ and let $\{v_1,\ldots,v_m\}$ be the rows of the matrix $M$. Then, consider the following matrices $\left(N = \left\lceil\frac{m}{r-1}\right\rceil\right)$:
   $$A_1 = \begin{pmatrix} v_1 \\ \vdots \\ v_{r-1} \\ {\sf 1}_{n-r+1,n} \end{pmatrix},~A_2 = \begin{pmatrix} {\sf 1}_{r-1,n} \\ v_r \\ \vdots \\ v_{2r-1} \\ {\sf 1}_{n-2r+1,n} \end{pmatrix}, \ldots , ~A_{N} = \begin{pmatrix} {\sf 1}_{(N-1)r-1,n} \\ v_{(N-1)r} \\ \vdots \\ v_{n}\end{pmatrix}.$$Then, it is easy to check that $M = A_1\star\cdots\star A_N$. 
    
    If $n \leq m$, we do the same constructions, considering columns instead of rows. 
   \qed
 \end{proof}
 Therefore, it makes sense to give the following definitions.
 \begin{definition}
  Let $M$ be a matrix and fix $r \geq 2$. We call an {\it $r$-th Hadamard decomposition} of $M$ an expression of the type $   M = A_1\star\ldots\star A_s,\text{ where } {\rm rk}(A_i) \leq r.$ We define the {\it $r$-th Hadamard rank} of $M$ as the smallest length of such a decomposition, i.e.,
  $$
   \Hrk_r(M) = \min\{ s ~|~ \text{there exist}~A_1,\ldots,A_s,~{\rm rk}(A_i) \leq r,~ M = A_1\star\ldots\star A_s\}.
  $$
%  Let $X_r \subset \PP({\rm Mat}_{m,n})$ be the space of matrices of rank at most $r$. We define the {\it $r$-th Hadamard generic rank} of $M$ as the smallest $s$ such that $M$ lies on the $s$-th Hadamard power of $X_r$, i.e.,
%  $$
%   \overline{\Hrk}_r(M) = \min\{ s ~|~ M \in X_r^{\star s}\}.
%  $$
  We define the {\it generic $r$-th Hadamard rank} of matrices of size $m\times n$ as
  $$
   \Hrk_r^\circ(m,n) = \min\{s ~|~ X_r^{\star s} = \PP({\rm Mat}_{m,n})\},
  $$
  and the {\it maximal $r$-th Hadamard rank} of matrices of size $m\times n$ as
  $$
   \Hrk_r^{\rm max}(m,n) = \max\{\Hrk_r(M) ~|~ M \in {\rm Mat}_{m,n}\}.
  $$
 \end{definition}
 We remark that these definitions may be seen as the multiplicative versions of the more common notion of {\it tensor ranks}, where we consider {\it additive decompositions} of tensors as sums of decomposable tensors. In terms of matrices, we look at decomposition as sums of rank $1$ matrices. A massive amount of work has been devoted to problems related to tensor ranks during the last few decades, especially due to their applications to statistics, data analysis, signal process, and others. See \cite{L} for a complete exposition of the current state of the art.

Hadamard product of matrices, i.e., the entrywise product, is the na\"ive definition for matrix multiplication that any school student would hope to study. Even if it is not the standard multiplication we have been taught, it is a very interesting operation, with nice properties and applications in matrix analysis, statistics and physiscs. As mentioned in the introduction, the generalization to the case of tensors has been used in data mining and quantum information \cite{CMS,PHYS}. We look at it from a geometric point of view, by studying Hadamard powers of varieties of matrices.

\smallskip

For a fixed positive integer $r \leq \min\{m,n\}$, we denote by $X_r \subset \PP({\rm Mat}_{m,n})$ the variety of matrices of size $m\times n$ with rank at most $r$. In other words, $X_r$ is the $r$-th secant variety of the Segre product $\PP^{m-1}\times \PP^{n-1}$. These are well-studied classic objects. Since ${\sf 1}_{m,n}$, the matrix of all $1$'s, which is the identity element for the Hadamard product, is contained in the variety $X_r$, we have a chain of inclusions as in \eqref{eq:chain}. 

% \begin{question}
%  Can we cover the whole space of matrices $\PP({\rm Mat}_{m,n})$ by considering a sufficiently high Hadamard power of $X_r$?
% \end{question}
 
 \begin{remark}
  Our aim is to study Hadamard powers of the varieties $X_r$ of matrices with rank at most $r$. As we observed before, we can view the Hadamard power $X_r^{\star 2}$ as a linear projection of the Segre product $X_r \times X_r$.  In terms of matrices, this is the geometric translation of the well-known fact that {\it the Hadamard product of two matrices is a submatrix of their Kronecker product}. Indeed, if $M=(m_{i,j})\in{\rm Mat}_{m,n}$ and $N=(n_{i,j})\in{\rm Mat}_{m,n}$, we define the Kronecker product as $M\otimes N = (m_{i,j}n_{h,k}) \in {\rm Mat}_{m^2,n^2}$. Then, $M\star N = (M\otimes N)|_{I,J}$, where $(M\otimes N)|_{I,J}$ denotes the restriction on the indexes $I = \{1,m+2,2m+3,\ldots,m^2\}$ and $J = \{1,n+2,2n+3,\ldots,n^2\}$.
 \end{remark}

Hadamard powers of a specific space of tensors has been considered in \cite{CMS} as the geometric interpretation of a particular statistical model. Therefore, we believe that the definitions of Hadamard ranks of matrices, and more generally of tensors, are very natural and may be an interesting area of research from several perspectives.

\smallskip
Proposition \ref{lemma:max rank} gives us an upper bound on the $r$-th Hadamard rank, i.e.,
\begin{center}$
 \Hrk_r^{\max}(m,n) \leq \left\lceil \frac{\min\{m,n\}}{r-1}\right\rceil.
$\end{center}
We can also give a lower bound on the generic rank as a straightforward application of the following well-known property of Hadamard product of matrices.
\begin{lemma}\label{lemma:bound_rank_Hadamard}
Given two matrices $A,B$, we have that 
  \begin{center}$
   {\rm rk}(A\star B) \leq {\rm rk}(A){\rm rk}(B).
  $\end{center}
\end{lemma}
\begin{proof}
  Say that ${\rm rk}(A) = r_1$ and ${\rm rk}(B) = r_2$. Consider the additive decomposition of $A$ and $B$ as sums of rank $1$ matrices, i.e.,
  \begin{center}$
   A = \sum_{i=1}^{r_1} a_i\cdot b_i^{\rm T} \text{ and } 
   B = \sum_{j=1}^{r_2} c_j\cdot d_j^{\rm T},
  $\end{center}
  where $a_i,b_i,c_j,d_j$ are column vectors. Then, we get that
  \begin{center}$
   A\star B = \sum_{i = 1}^{r_1} \sum_{j=1}^{r_2} (a_i \star c_j)\cdot (b_i\star d_j)^{\rm T}.
  $\end{center}
  Therefore, we have that ${\rm rk}(A\star B) \leq r_1r_2$.
  \qed
 \end{proof}
% This lemma can be rephrased by saying that $X_r^{\star 2} \subset X_{r^2}$,  for any $r$.
As an immediate consequence of this lemma we see that that $X_r^{\star 2} \subset X_{r^2}$, for any $r$.
 In particular, we obtain a lower bound on the generic Hadamard rank.
  \begin{corollary}
  Fix $r\geq 2$. Then, the generic $r$-th Hadamard rank of matrices of size $m \times n$ is at least $\left\lceil\log_r(\min\{m,n\})\right\rceil$.
 \end{corollary}
 \begin{proof}
  If $s <  \left\lceil\log_r(\min\{m,n\})\right\rceil$, then $r^s < \min\{m,n\}$. Hence, the Hadamard product of $s$ matrices of rank $r$ cannot have maximal rank and, therefore, it cannot be enough to cover all the space of matrices of size $m\times n$.
 \end{proof}
Therefore, we have the following chain of inequalities.
\begin{equation}\label{ineq:trivial}
\left\lceil\log_r(\min\{m,n\})\right\rceil \leq \Hrk_r^{\circ}(m,n) \leq \Hrk^{\max}_r(m,n) \leq 
\left\lceil \frac{\min\{m,n\}}{r-1}\right\rceil.
\end{equation}
By this chain of inclusions we get the following result.
\begin{proposition}\label{prop:m-1_Hrk}
 Let $m \leq n$ and consider $r = m - 1$. Then, we have 
 \begin{center}$
  \Hrk^{\circ}_{m-1}(m,n) = \Hrk^{\max}_{m-1}(m,n) = 2.
 $\end{center}
\end{proposition}
\begin{proof}
 On the left hand side of \eqref{ineq:trivial} we have $\left\lceil\log_{m-1}(m)\right\rceil = 2.$
 
 On the right hand side, we have $\left\lceil \frac{m}{m-2}\right\rceil$, which is equal to $2$ if $m\geq 4$. Then, in order to conclude, we just need to prove the case $m = 3$.
 
 Let $m = 3$. If we consider a matrix $M$ of rank $\leq 2$, then it lies on $X_2$. Assume that $M$ has rank $3$ and let $v_i = (v_{i,1},\ldots,v_{i,n})$, for $i = 1,2,3$, be the rows of $M$.  Consider the first two rows. If $v_{1,j}$ and $v_{2,j}$ are not both equal to zero, for all $j = 1,\ldots,n$, then there exists a linear combination of $\lambda v_1 + \mu v_2$ with all entries different from zero and, therefore, we can decompose $M$ as follows
 $$
  M = \begin{pmatrix} v_{1,1} & \ldots & v_{1,n} \\ v_{2,1} & \ldots & v_{2,n} \\ \lambda v_{1,1} + \mu v_{2,1} & \ldots & \lambda v_{1,n} + \mu v_{2,n}\end{pmatrix} \star \begin{pmatrix} 1 & \ldots & 1 \\ 1 & \ldots & 1 \\ \frac{v_{3,1}}{\lambda v_{1,1} + \mu v_{2,1}} & \ldots & \frac{v_{3,n}}{\lambda v_{1,n} + \mu v_{2,n}}\end{pmatrix}.
 $$
 If we have $v_{1,j} = v_{2,j} = 0$, for some $j = 1,\ldots,n$, any linear combination of $v_1$ and $v_2$ will have the $j$-th entry equal to zero. Therefore, we cannot use the previous algorithm. Hence, we define $\widetilde{v}_i$, for $i = 1,2$, as
 \begin{center}$
  \widetilde{v}_{i,j} = \begin{cases} v_{i,j} & \text{ if } v_{1,j} \neq 0 \text{ or } v_{2,j} \neq 0; \\
  1 & \text{ if } v_{1,j} = v_{2,j} = 0. \end{cases}
 $\end{center}
 Now, there exists a linear combination of $\lambda \widetilde{v}_1 + \mu \widetilde{v}_2$ with all entries different from zero. Therefore, if we define a row $u$ as
 \begin{center}$
  u_i = \begin{cases} 1 & \text{ if } v_{1,j} \neq 0 \text{ or } v_{2,j} \neq 0; \\
  0 & \text{ if } v_{1,j} = v_{2,j} = 0,\end{cases}
 $\end{center}
 we can decompose $M$ as 
  \begin{center}$
  M = \begin{pmatrix} \widetilde{v}_{1,1} & \ldots & \widetilde{v}_{1,n} \\ \widetilde{v}_{2,1} & \ldots & \widetilde{v}_{2,n} \\ \lambda \widetilde{v}_{1,1} + \mu \widetilde{v}_{2,1} & \ldots & \lambda \widetilde{v}_{1,n} + \mu \widetilde{v}_{2,n}\end{pmatrix} \star \begin{pmatrix} u_1 & \ldots & u_n \\ u_1 & \ldots & u_n \\ \frac{v_{3,1}}{\lambda \widetilde{v}_{1,1} + \mu \widetilde{v}_{2,1}} & \ldots & \frac{v_{3,n}}{\lambda \widetilde{v}_{1,n} + \mu \widetilde{v}_{2,n}}\end{pmatrix}.
 $\end{center}
 Therefore, $\Hrk_2^{\max}(3,n) = 2$.
 \qed
 \end{proof}
 \begin{example}
  Consider the matrix $M = \begin{pmatrix}1 & 2 & 0 & 1\\ -1 & 1 & 0 & 0\\ 0 & 1 & 1 & 2 \end{pmatrix}$. Then, we consider
  \begin{center}$
   \widetilde{v}_1 = (1,2,1,1),~\widetilde{v}_2 = (-1,1,1,0),~u=(1,1,0,1).
  $\end{center}
  Hence,
  \begin{center}$
   M = \begin{pmatrix}1 & 2 & 1 & 1 \\ -1 & 1 & 1 & 0 \\ 1 & 5 & 3 & 2\end{pmatrix} \star \begin{pmatrix}1 & 1 & 0 & 1\\ 1 & 1 & 0 & 1\\ 0 & \frac{1}{5} & \frac{1}{3} & 1 \end{pmatrix}.
  $\end{center}
 \end{example}
 \begin{remark}\label{rmk:upper_bound}
 We proved that for $r = \min\{m,n\}-1$, the $r$-th Hadamard rank is equal to $2$. Actually, the upper-bound in \eqref{ineq:trivial} let us be more precise. Indeed, we can say that for any $\frac{\min\{m,n\}+2}{2} < r < \min\{m,n\}$, we get $\Hrk_r^{\circ}(m,n)=2$.
 \end{remark}
 
 In other cases, we need a more geometric approach in order to understand the generic Hadamard rank. 
By using Proposition \ref{prop:expdim}, we can define the expected dimension for the $s$-th Hadamard power of the variety $X_r$ of rank $r$ matrices.
\begin{proposition}
 In the same above notation,
 \begin{equation}\label{eq:exp_dim}
  \dim(X_r^{\star s}) \leq \min\big\{s\dim(X_r) - (s-1)\dim(X_1),\dim~\PP({\rm Mat_{m,n}})\big\}.
 \end{equation}
\end{proposition}
\begin{proof}
We proceed by induction on $s$. For $s = 1$, it follows trivially from definitions. Consider $s > 1$. Then, since $X_r^{\star s} = X_r^{\star(s-1)}\star X_r$, by Proposition \ref{prop:expdim} and by inductive hypothesis, we get
\begin{align*}
  \dim(X_r^{\star s}) & \leq \min\big\{\dim(X_r^{\star (s-1)}) + \dim(X_r) - \dim(X_1),\dim\PP({\rm Mat_{m,n}})\big\} \\
  & = \min\big\{s\dim(X_r) - (s-1)\dim(X_1),\dim\PP({\rm Mat_{m,n}})\big\}.
\end{align*}
\qed
\end{proof}
We refer to the formula on the right hand side of \eqref{eq:exp_dim} as the {\it expected dimension} of $X_r^{\star s}$. More precisely, we have the following
\begin{align*}
 {\rm exp}.\dim(X_r^{\star s}) & = \min\big\{s\dim(X_r) - (s-1)\dim(X_1),\dim\PP({\rm Mat_{m,n}})\big\}  \\
  & = \min\big\{ sr(n+m-r) - (s-1)(n+m-1),mn\big\}-1.
\end{align*}Therefore, the {\it expected generic $r$-th Hadamard rank} is
\begin{equation}\label{eq:expected_Hadamard_rank}
 {\rm exp}.\Hrk^\circ_r(m,n) = \left\lceil\frac{\dim\PP({\rm Mat}_{m,n}) - \dim(X_1)}{\dim(X_r) - \dim(X_1)}\right\rceil = \left\lceil\frac{mn-(m+n-1)}{r(m+n-r)-m-n+1}\right\rceil  
\end{equation}
\begin{remark}
 A very important concept in the world of tensors additive decomposition is the idea of {\it identifiability}, namely, we say that a tensor is {\it identifiable} if it has a unique decomposition as sum of decomposable tensors. Since we are viewing Hadamard decomposition as a multiplicative version of tensor decomposition, we might look for identifiability also in this set up. However, in this case, we cannot have identifiability for any matrix. Indeed, consider a $r$-th Hadamard decomposition of a matrix $M$, i.e., we have
 \begin{center}$
   M = A_1 \star \cdots \star A_s,\text{ with }{\rm rk}(A_i) = r;
 $\end{center}
 then, for any $(s-1)$-tuple of rank $1$ matrices $R_1,\ldots,R_{s-1}$, all with non-zero entries, we can construct a different $r$-th Hadamard decomposition as
 \begin{center}$
   M = \big(R_1\star A_1\big) \star \cdots \star \big(R_{s-1}\star A_{s-1}\big) \star \big((R_1\star\cdots\star R_{s-1})^{\star (-1)}\star A_s\big),
 $\end{center}
 where $R^{\star (-1)}$ denotes the Hadamard inverse of the matrix $R$. Here, we have to recall that ${\rm rk}(R_i\star A_i) \leq {\rm rk}(A_i)$, for any $i = 1,\ldots,s-1$, by Lemma \ref{lemma:bound_rank_Hadamard}, and, similarly, we have ${\rm rk}\big((R_1\star\cdots\star R_{s-1})^{\star (-1)}\star A_s\big) \leq {\rm rk}(A_s)$, because ${\rm rk}(R_1\star\cdots\star R_{s-1})^{\star (-1)} = 1$.
\end{remark}
We can check that \eqref{eq:exp_dim} is the actual dimension and, consequently, \eqref{eq:expected_Hadamard_rank} gives the correct generic $r$-th Hadamard rank for matrices of small size.

Here we describe an algorithm written with \emph{Macaulay2} to compute the dimensions of Hadamard powers of varieties of square matrices of given rank. This allows us to compute the corresponding generic Hadamard ranks (Table \ref{table}). We reduced to square matrices for simplicity of exposition, but the code can be easily generalized.

The key point is to use Lemma \ref{lemma:Terracini} which states that the tangent space to $X_r^{\star s}$ at a generic point $A_1\star \cdots \star A_s$ is given by
 \begin{equation}\label{formula:tangent_space}
  T_{A_1\star \cdots \star A_s}(X_r^{\star s}) = \left\langle T_{A_1}(X_r)\star A_2\star \cdots \star A_{s},\ldots,A_1\star \cdots \star A_{s-1} \star T_{A_s}(X_r)\right\rangle
 \end{equation}

Hence, we first need to construct the tangent spaces at $s$ random points of $X_r$.

 Recall that, if $A$ is a matrix of rank $r$ written as $A = \sum_{i=1}^r u_i\cdot v_i^T$, $u_i,v_i\in\CC^n$, the tangent space of $X_r$ at $A$ is given by  
 \begin{center}
 $
  T_A\left(X_r\right) = \left\langle u_1\cdot(\CC^n)^T + (\CC^n)\cdot v_1^T,\ldots,u_r\cdot(\CC^n)^T + (\CC^n)\cdot v_r^T\right\rangle.
 $
 \end{center}
 \newpage
 \noindent Here is the {\it Macaulay2} code.
{\small 
 \begin{verbatim}
-----------------------------------------------------------
INPUT:  n = sizes of matrices;
        r = rank of matrices;
        s = Hadamard power to compute;
        
OUTPUT: D = dimension of the s-th Hadamard power of
           the variety of rank r matrices of size nxn.
-----------------------------------------------------------
   S := QQ[z_(1,1)..z_(n,n),
            a_(1,1)..a_(n,r),b_(1,1)..b_(n,r),
            c_(1,1)..c_(2*r,n)];
---- Construct s random matrices of rank r
    u = for i from 1 to s list 
            for j from 1 to 2*r list 
                random(S^n,S^{0});
    A = for i from 0 to (s-1) list sum (
             for j from 0 to (r-1) list 
                 u_i_(2*j) * transpose(u_i_(2*j+1))
        );
---- Construct their tangent spaces
    C = for i from 1 to 2*r list 
            genericMatrix(S,c_(i,1),n,1);
    TA = for i from 0 to (s-1) list sum 
             for j from 0 to (r-1) list 
                 u_i_(2*j) * transpose C_(2*j) + 
                     C_(2*j+1) * transpose(u_i_(2*j+1));
 \end{verbatim}
 }
 Now, we construct the vector spaces spanning the tangent space of $X_r^{\star s}$ as in \eqref{formula:tangent_space}. First, we define a function \texttt{HP} to compute the Hadamard product of two matrices.
 {\small
 \begin{verbatim}
---- Method to construct the Hadamard product of a
---- list of matrices of same size;
   HP = method();
   HP List := L -> (
    s := #L;
    r := numRows(L_0);
    c := numColumns(L_0);
    for i from 1 to (s-1) do 
      if (numRows(L_i)!=r or numColumns(L_i)!=c) then 
          return << "error";
    H := for i from 0 to (r-1) list 
           for j from 0 to (c-1) list product (
              for h from 0 to (s-1) list (L_h)_j_i
              );         
    return matrix H
    )
---- Construct the two vector spaces spanning the tangent  
---- space of the Hadamard power and find their equations 
---- in the space of matrices
   TAstar = for i from 0 to (s-1) list 
                HP(toList(set{TA_i}+set(A)-set{A_i}));
   M = genericMatrix(S,z_(1,1),n,n);
   H = for i from 0 to (s-1) list 
            ideal flatten entries (M - TAstar_i);
   H1 = for i from 0 to (s-1) list 
             eliminate(toList(c_(1,1)..c_(2*r,n)),H_i);
   T = QQ[z_(1,1)..z_(n,n)];
   E = for i from 0 to (s-1) list sub(H1_i,T);
 \end{verbatim}
 }
 In \texttt{E}, we have the list of the equations of the tangent spaces to the variety $X_r$ at the $s$ random points. From these, we can construct a vector basis for each tangent space. Now, in order to compute the dimension of their span it is enough to compute the rank of the matrix obtained by collecting all these vector basis together.
 {\small 
 \begin{verbatim}
 K = for i from 0 to (s-1) list 
        kernel transpose 
           contract(transpose vars(T),mingens E_i);
 tt = mingens K_0 | mingens K_1;
 if s >= 3 then (
        for i from 2 to (s-1) do tt = tt | mingens K_i
        );
 D = rank tt
 \end{verbatim}
 }
\noindent  In the following table, we list the generic $r$-th Hadamard ranks that we have computed for square matrices of small size. 
 \begin{table}
 \begin{center}
 \begin{tabular}{c|c|c}
   $~n~$ & $~r~$ & $r$-th Hadamard rank \\
   \hline\hline
   3 & 2 & 2 \\
   4 & 2 & 2 \\
   5 & 2 & 3 \\
   & 3 & 2 \\
   6 & 2 & 3 \\
   & 3 & 2 \\
   7 & 2 & 4 \\
   & 3 & 2 \\
   & 4 & 2 \\
   8 & 2 & 4 \\
   & 3 & 3 \\
   & 4 & 2 \\
   9 & 2 & 5 \\
   & 3 & 3 \\
  \end{tabular}~~~~~~
   \begin{tabular}{c|c|c}
   $~n~$ & $~r~$ & $r$-th Hadamard rank \\
   \hline\hline
   9 & 4 & 2 \\
   & 5 & 2 \\
   10 & 2 & 5 \\
   & 3 & 3 \\
   & 4 & 2 \\
   & 5 & 2 \\
   11 & 2 & 6 \\
   & 3 & 3 \\
   & 4 & 2 \\
   & 5 & 2 \\
   & 6 & 2 \\
   12 & 2 & 6 \\
   & 3 & 4 \\
   & 4 & 3 \\
  \end{tabular}~~~~~~\begin{tabular}{c|c|c}
   $~n~$ & $~r~$ & $r$-th Hadamard rank \\
   \hline\hline
  12  & 5 & 2 \\ 
   & 6 & 2 \\
   13& 2 & 7 \\
   & 3 & 4 \\
   & 4 & 3 \\
   & 5 & 2 \\
   & 6 & 2 \\
   & 7 & 2 \\
   14 & 2 & 7 \\
   & 3 & 4 \\
   & 4 & 3 \\
   & 5 & 2 \\
   & 6 & 2 \\
   & 7 & 2 \\
  \end{tabular}
  \caption{Generic $r$-th Hadamard ranks of square matrices of size $n\times n$ with $n \leq 14$. By Remark \ref{rmk:upper_bound}, we could restrict to the cases $r < \frac{n+2}{2}$; for $r \geq \frac{n+2}{2}$, we know that $\Hrk_r^\circ(n,n) = 2$. This computation required less than 9 minutes on a laptop with a processor 2,2GHs Intel Core i7.}\label{table}
 \end{center} 
 \end{table}
\begin{acknowledgement}
 This article was initiated during the Apprenticeship Weeks (22 August-2 September 2016), led by Bernd Sturmfels, as part of the Combinatorial Algebraic Geometry Semester at the Fields Institute. The second author was supported by {\it G S Magnuson Foundation} from {\it Kungliga Vetenskapsakademien} (Sweden).
\end{acknowledgement}

\end{document}